\documentclass[11pt]{amsart}
\usepackage{geometry}                
\usepackage{graphicx}
\usepackage{amssymb}
\usepackage{epstopdf}
 \usepackage{latexsym}
\usepackage{amsmath}
\usepackage{amsfonts} 
\usepackage{amsthm}

\newtheorem{Proposition}{Proposition}

    \def\sqr#1#2{{\vcenter{\vbox{\hrule height .#2pt
                             \hbox{\vrule width .#2pt height#1pt \kern#1pt
                                   \vrule width .#2pt}
                             \hrule height .#2pt}}}}

     \def\CC{\mathbb{C}}
    
    \def\NN{\mathbb{N}}

    \def\ZZ{\mathbb{Z}}

 \def\bfu{\mathbf{u}}
 \def\bfw{\mathbf{w}}
 \def\bff{\mathbf{f}}
  \def\bfg{\mathbf{g}}
 
 \def\bfy{\mathbf{y}}
  \def\bfh{\mathbf{h}}
    \def\bfP{\mathbf{P}}
 
  \def\bfq{\mathbf{q}}
   \def\bfm{\mathbf{m}}   
      
\def\bflam{\boldsymbol{\lambda}}

\def\be{\begin{equation}}
\def\ee{\end{equation}}

\begin{document}

\title[Generalization of Jacobi Polynomials]{A class of matrix-valued polynomials generalizing Jacobi Polynomials}
\author{Rodica D. Costin}


\maketitle

\ 

\

\begin{abstract}

A hierarchy of matrix-valued polynomials which generalize the Jacobi polynomials is found.  
Defined by a Rodrigues formula, they are also products of a sequence of differential operators. Each class of polynomials is complete, satisfies a two-step recurrence relation, integral inter-relations, and quasi-orthogonality relations.

\end{abstract}

\

\section{Motivation}\label{Introduction}

The understanding of matrix-valued orthogonal polynomials has advanced greatly in recent years, and a wealth of references can be found in the recent work of Barry Simon \cite{Simon}.

The polynomials introduced in the present paper generalize the Jacobi polynomials.
A different generalization, yielding an orthogonal class, is due to Gr\"unbaum \cite{Grunbaum}.

The polynomials defined here were not sought-for. Rather they appear naturally in the forefront of a different problem, the study of the possibility to linearize (by an analytic change of variables) a differential equation whose linear part has only regular singular points.  In a neighborhood of one such point equations are (generically) linearizable \cite{RDC_Euler}. But this is no longer the case if the domain studied contains two such singularities. For example:
\be\label{eq}
\frac{d\bfu}{dx}=\, M\, \bfu\,  +\, \bff(x,\bfu)\ \ \ \ {\mbox{with}}\ M=\frac{1}{x-1}\, A+\, \frac{1}{x+1}\, B
\ee
($\bff$ collects all the nonlinear terms in $\bfu$) is {\em{not}} necessarily equivalent to its linear part
\be\label{lineq}
\frac{d\bfw}{dx}=\, M\, \bfw
\ee
for $x$ in a domain in the complex plane including both singular points $\pm 1$. 

It turns out, however that for any nonlinear term $\bff(x,\bfu)$ there exists a unique $\phi(\bfu)$ so that the equation with the "corrected" nonlinear part $\bff(x,\bfu)-\phi(\bfu)$ is linearizable \cite{Norm_Form}.

Besides its clear intrinsic interest, the problem of detecting linearizable equations is important also since linearizability and integrability turn out to be intimately connected \cite{RDC_MDK}.

Looking for an analytic change of variables $\bfy=\bfw+\bfh (x,\bfw)$ that transforms (\ref{lineq}) into 
(\ref{eq})  one is lead to the study of the homological equation which boils down to equations of the form
\be\label{homeq}
\partial_x\bfh+d_\bfw\bfh\, M\bfw-M\bfh=\bfg(x,\bfw)/(x^2-1)
\ee
where now all functions can be assumed homogeneous polynomials in $\bfw$, of degree, say, $n$. 

The main question is: for which functions $\bfg$ does equation (\ref{homeq}) possess a solution 
$\bfh$ which is analytic in a domain containing both singularities $x=\pm 1$?

In the scalar case, when $\bfw=w\in\CC$, $\bfh=h, \bfg=g\in\CC$, and $A=a,B=b\in\CC$, it turns out that: equation  (\ref{homeq}) has an analytic solution if and only if the first coefficient in the Jacobi series in $\left\{ P_k^{((n-1)a-1,(n-1)b-1)}\right\}_{k=0,1,2,\ldots}$ of $g$ is zero. Moreover, Jacobi series expansions are canonical to approaching the problem, since for $g(x)=P_{k+1}^{((n-1)a-1,(n-1)b-1)}(x)$ the unique analytic solution of  (\ref{homeq}) is $h=\frac{-1}{2k}P_{k}^{((n-1)a,(n-1)b)}$.

In the linearization problem mentioned the double sequence of Jacobi polynomials $P_{k}^{(na,nb)}$ for $k\geq 0,n\geq 2$ is interconnected. Numerical expansions in Jacobi series yield (what appear to be) rapidly convergent series.\footnote{However,  the proof of convergence of the series obtained, which did not use the Jacobi structure, turned out to be extremely delicate: a steepest descent method was used  and small denominators had to be dealt with \cite{RDC_AnLin}.}

The natural question was whether the polynomial structure found in the one-dimensional case survives in more dimensions. 

The answer turns out affirmative - see (\ref{inteq1}), (\ref{inteq2}). These multidimensional polynomials share many properties common to the usual classes of orthogonal polynomials: each class is complete, satisfies a two-step recurrence relation, quasi-orthogonality relations exist, and  
there are integral inter-relations. It is the author's belief that more that what is presented here does hold.

\section{Main Results}

\subsection{Notations}\label{notat}

Consider a Fuchsian differential equation in $\CC^d$ with three singularities in the extended complex domain $\overline{\CC}$:
\be\label{Fuchs_eq}
\bfy '\, =M\bfy\ \ \ \ \ {\mbox{where}}\ \ \ \  M=\frac{1}{x-1}\, A\, +\, \frac{1}{x+1}\, B\ ,\ \ (x\in\CC,\ \bfy\in\CC^d)
\ee
where $A$ and $B$ are $d\times d$ matrices so that \newline
(i) $A+B$ is diagonalizable, and \newline
(ii) $A+B$ satisfies the following
{\bf{ nonresonance condition:}} the eigenvalues $\lambda_1,\ldots ,\lambda_d$ of the matrix $A+B$ satisfy:
\be\label{non_res}
\mathbf{n}\cdot {\bflam}-\lambda_j \not\in\ZZ_-\ \ \ {\mbox{for\ all}}\ \ \mathbf{n}\in\NN^d\ ,\ \ j=1,\ldots ,d
\ee 
These assumptions are, of course, satisfied by generic matrices.

Denote for simplicity:
\be\label{nota}
M_1=A+B,\ \ \ M_2=A-B,\ \ \ Q=x^2-1
\ee
so that $ QM=xM_1+M_2$.

It will be assumed, without loss of generality, that the matrix $M_1$ is diagonal.

Let $Y$ be a fundamental matrix of solutions of (\ref{Fuchs_eq}), therefore 
$$Y'=MY\ {\mbox{and}}\ \left(Y^{-1}\right)'=-Y^{-1}M$$

Denote by $\mathcal{P}_n$ the space of vector-valued polynomials in $\bfw=(w_1,\ldots,w_d)$, homogeneous, degree n:
$$ \mathcal{P}_n\, =\, \left\{\ \bfq\ ;\ \bfq(\bfw)=\sum_{\bfm\in\NN^d, |\bfm|=n}\bfq_\bfm\bfw^\bfm,\ \bfq_\bfm\in\CC^d\ \right\}$$
where $|\bfm|=m_1+\ldots+m_d$, $\bfw^\bfm=w_1^{m_1}\ldots w_d^{m_d}$, $\NN=\{0,1,2,\ldots\}$.

\subsection{The class of polynomials: definition}\label{definition}

The polynomial in $x$ degree k, $\bfP_k(x)$,
is, for every $x$, the linear operator on $\mathcal{P}_n$ defined by
\be\label{defPk}
\bfP_k(x)\bfq(\bfw)\, =\, Y(x)\, \frac{d^k}{dt^k}\Big| _{t=x}\, \left[\, Q(t)^k\, Y^{-1}(t)\, \bfq\left(Y(t)Y^{-1}(x)\bfw\right)\,\right]
\ee

Formula (\ref{defPk})  can be compactly written as follows.  Denoting by $W(x)$ the linear operator on $\mathcal{P}_n$:
$$W(x)\bfq(\bfw)\, =\, Y(x)^{-1}\, \bfq\left(Y(x)\bfw\right)$$
(which gives $\bfq$ in new, $x$-dependent coordinates in $\CC^d$) the definition (\ref{defPk}) is
\be\label{defPkW}
\bfP_k(x)\, =\, W(x)^{-1}\frac{d^k}{dx^k}\left[ Q(x)^kW(x)\right]
\ee
which is a Rodrigues formula. 

The fact that $\bfP_k(x)$ are polynomials in $x$, degree $k$, is shown in 
Proposition\,\ref{Ppoly}.

Denote by $D_1,D_2$ the following linear operators on $\mathcal{P}_n$:
\be\label{defD12}
D_i\bfq(\bfw)\, =\, {\rm{d}}\bfq(\bfw)\, M_i\bfw-M_i\bfq(\bfw),\ \ i=1,2
\ee
where $M_i$ are defined in (\ref{nota}).

\begin{Proposition}\label{Ppoly}

The operators $\bfP_k(x)$ are polynomials in $x$ degree k. 

Moreover, they can be written as the composition 
\be\label{prodPk}
\bfP_k=\mathcal{A}_1\mathcal{A}_2\ldots\mathcal{A}_k
\ee
where\footnote{Note that the action of first operator $\mathcal{A}_k$ when applied to $\bfq(\bfw)$ (which does not depend on $x$) equals $[x\left(2j+D_1\right)+D_2]\bfq(\bfw)$.}
\be\label{defAj}
\mathcal{A}_j=x\left(2j+D_1\right)+D_2+Q\partial_x
\ee

\end{Proposition}

{\em{Proof.}}

Note that calculating the first derivative in (\ref{defPk}) we obtain
$$\frac{d}{dt}\left[Q(t)^kY(t)^{-1}\bfq(Y(t)\bfu)\right]=Q(t)^{k-1}Y(t)^{-1}\left(2kt+tD_1+D_2\right)\bfq(Y(t)\bfu)$$
$$\equiv Q(t)^{k-1}Y(t)^{-1}\mathbf{r}(t,Y(t)\bfu)\ \ \ {\mbox{where\ }}\mathbf{r}(t,\bfw)=\mathcal{A}_k\bfq(\bfw)$$

Then
$$\frac{d^2}{dt^2}\left[Q(t)^kY(t)^{-1}\bfq(Y(t)\bfu)\right]=\frac{d}{dt}\left[Q(t)^{k-1}Y(t)^{-1}\mathbf{r}(t,Y(t)\bfu)\right]$$
$$=Q(t)^{k-2}Y(t)^{-1}\left[2(k-1)t+tD_1+D_2+Q\partial_t\right]\mathbf{r}(t,Y(t)\bfu)$$
$$=Q(t)^{k-2}Y(t)^{-1}\mathcal{A}_{k-1}\mathbf{r}\left( t,Y(t)\bfu\right)\, =\, Q(t)^{k-2}\,Y(t)^{-1}\ \mathcal{A}_{k-1}\mathcal{A}_k\bfq\left( Y(t)\bfu\right)$$

The $k$ derivatives in (\ref{defPk}) can be calculated in this way recursively, yielding (\ref{prodPk}) and completing the proof of Proposition \ref{Ppoly}. \qed

{\bf{Remark 1.}} The same definition (\ref{defPk}) can be used to define $\bfP_k(x)$ as linear operators on the space of $\CC^d$-valued formal series in $\bfw$, or on the space of convergent such series.

\subsection{The polynomials (\ref{defPk}) generalize the Jacobi polynomials}

Consider the one-dimensional case: $d=1$. Equation (\ref{Fuchs_eq}) and a fundamental solution are, in this scalar case,
$$y'=\left(\frac{a}{x-1}+\frac{b}{x+1}\right)y\ \ \ \ \ {\mbox{and\ \ \ \ }}\ y(x)=(x-1)^a(x+1)^b$$
The homogeneous polynomials degree $n$ are just multiples of $w^n$: $q(w)=c\, w^n$ and formula 
(\ref{defPk}) gives
$$P_k(x)q=\frac{1}{y(x)^{n-1}}\, \frac{d^k}{dx^k}\left[Q(x)^ky(x)^{n-1}\right]\, q$$
which is the Rodrigues formula for the Jacobi polynomial $P_k^{((n-1)a,(n-1)b)}(x)$ (up to a multiplicative factor).

{\bf{Remark 2.}} In the one-dimensional case the operators $D_{1,2}$ of (\ref{defD12}) are multiplication by $(n-1)(a+b)$, and $(n-1)(a-b)$ respectively, and formula (\ref{prodPk}) gives the following representation for Jacobi polynomials (up to a numerical factor):
\be\label{prodJacobi}
P_k^{(\alpha,\beta)}(x)=\mathcal{A}_1\mathcal{A}_2\ldots\mathcal{A}_k\, 1
\ee
where $\mathcal{A}_j=(2j+\alpha+\beta)x+(\alpha-\beta)+Q\partial_x$

Formula (\ref{prodJacobi}) provides the following elegant way to deduce the second order operator for which Jacobi polynomials are eigenvalues: noting that 
$$\partial_x\mathcal{A}_j=(\alpha+\beta+2j)+\mathcal{A}_{j+1}\partial_x$$
 we obtain
$$\mathcal{A}_1\partial_xP_k^{(\alpha,\beta)}=\mathcal{A}_1\partial_x\mathcal{A}_1\mathcal{A}_2\ldots\mathcal{A}_k=(\alpha+\beta+2)P_k^{(\alpha,\beta)}+\mathcal{A}_1\mathcal{A}_2\partial_x\mathcal{A}_3\ldots\mathcal{A}_k$$
$$=\ldots=\left[(\alpha+\beta+2)+\ldots+(\alpha+\beta+2k)\right]P_k^{(\alpha,\beta)}=k(\alpha+\beta+k+1)P_k^{(\alpha,\beta)}$$

\subsection{The multidimensional commutative case.} If the matrices $A$ and $B$ are simultaneously diagonalizable:
$$A\, =\, {\mbox{diag}}\, \left[a_i\right]_{i=1,\ldots , d}\ {\mbox{and}}\ B\, =\, {\mbox{diag}}\, \left[b_i\right]_{i=1,\ldots , d}$$
then the matrix $Y$ is
$$Y={\mbox{diag}}\, \left[\, (x-1)^{a_i}(x+1)^{b_i}\, \right]{i=1,\ldots , d}$$
and the action of the operators (\ref{defPk}) on the canonical basis  of $\mathcal{P}_n$: 
$\mathbf{b}_{\bfm,j}=\bfw^\bfm \mathbf{e}_j$, $j=1,\ldots ,d$, $|\bfm|=n$ is
\be\label{Jacobi}
\bfP_k(x)\mathbf{b}_{\bfm,j}\, =\, \frac{1}{W_{\bfm,j}(x)}\, \frac{d^k}{dx^k}\,\left[\, Q(x)^kW_{\bfm,j}(x)\, \right]\, \mathbf{b}_{\bfm,j}
\ee
where
\be\label{WJ}W_{\bfm,j}(x)=(x-1)^{\bfm\cdot\mathbf{a}-a_j}(x+1)^{\bfm\cdot\mathbf{b}-b_j}
\ee
therefore $\bfP_k(x)$ are diagonal operators, with entries Jacobi polynomials degree $k$.

This is a family of matrix-valued orthogonal polynomials with respect to the inner product
$$<\bfP_l,\bfP_k>\, :=\, \int_{-1}^1\, {\rm{Tr}}\, \left\{\, \overline{\bfP}_l(x)\, W(x)\, \bfP_k(x)\,\right\}\ dx$$
$$=\, \sum_{|\bfm|=n,j=1,\ldots,d}\, \int_{-1}^1\, \overline{P_l}^{(\bfm\cdot\mathbf{a}-a_j,\bfm\cdot\mathbf{b}-b_j)}{W_{\bfm,j}}\, {P_k}^{(\bfm\cdot\mathbf{a}-a_j,\bfm\cdot\mathbf{b}-b_j)}\ dx$$
in the cases when the integrals exist.\footnote{Appropriate conditions relating $\mathbf{a},\mathbf{b},n,d$ will ensure existence of the integrals.}

\subsection{The non-commutative linear case.} For $n=1$, the homogeneous polynomials degree one are linear transformations: $\bfq(\bfw)=q\bfw$ where $q$ is a matrix. Formula (\ref{defPk}) is, in this case,
\be\label{defPklin}
\bfP_k(x)q\, =\, Y(x)\, \frac{d^k}{dx^k}\, \left[\, Q(x)^k\, Y^{-1}(x)\, q\, Y(x)\,\right]\, Y^{-1}(x)
\ee
which implies
\be\label{lin_case}
{\rm{Tr}}\, \left[\, \bfP_k(x)q\,\right]\, =\, P_k^{(0,0)}(x)\, {\rm{Tr}}\, q
\ee
where $P_k^{(0,0)}(x)=P_k(x)$ are the Legendre polynomials.

\subsection{Properties of the generalized polynomials (\ref{defPk})}

\subsubsection{Rodrigues formula.} The definition (\ref{defPkW}) of $\bfP_k(x)$ is a Rodrigues formula. 

\subsubsection{A two-step recurrence relation}\label{two_step_rec}

\begin{Proposition}\label{recurrence}

The operators $\bfP_k(x)$ defined by  (\ref{defPk}) satisfy the following two-step recurrence relation
\be\label{recPk}
x\, \bfP_k(x)\, =\bfP_{k+1}(x)\alpha_k+\bfP_k(x)\, \beta_k+\bfP_{k-1}(x)\gamma_k
\ee
where $\alpha_k,\beta_k,\gamma_k$ are linear operators on $\mathcal{P}_n$, and in particular
\be\label{valal}
\alpha_k=(D_1+2k+1)^{-1}(D_1+2k+2)^{-1}(D_1+k+1)
\ee

\end{Proposition}

{\em{Proof of Proposition\,\ref{recurrence}.}}

Note the identities
\be\label{ide1}
\mathcal{A}_j(x\,\mathbf{r})\, =\, x\,\mathcal{A}_j(\mathbf{r})+Q\,\mathbf{r}\ \ \ \ \ \ \ \ \ \ \ \ \ \ \ \ \ \ \ \ \ \ 
\ee
and
\be\label{ide2}
Q\,\mathcal{A}_j(\mathbf{r})=\mathcal{A}_{j-1}(Q\,\mathbf{r})\ \ \ \ \ \  \ \ \ \ \  \ \ {\mbox{for\ any\ }} \mathbf{r}=\mathbf{r}(x,\cdot)\in\mathcal{P}_n\\
\ee
which by iteration give
$$
\mathcal{A}_1\mathcal{A}_2\ldots\mathcal{A}_k(x\,\bfq)=\mathcal{A}_1\mathcal{A}_2\ldots\mathcal{A}_{k-1}\left( x\,\mathcal{A}_k\bfq\right)\, +\,  \mathcal{A}_1\mathcal{A}_2\ldots\mathcal{A}_{k-1}\, (Q\,\bfq)$$
\be\label{ide3}
=\, \ldots\, =x\, \mathcal{A}_1\mathcal{A}_2\ldots\mathcal{A}_k\, \bfq\, +\, k\, \mathcal{A}_1\mathcal{A}_2\ldots\mathcal{A}_{k-1}\, (Q\,\bfq)
\ee
Using (\ref{prodPk}) and (\ref{ide3}) relation (\ref{recPk}) follows if we have
$$\mathcal{A}_{k}x-kQ=\mathcal{A}_k\mathcal{A}_{k+1}\alpha_k+\mathcal{A}_k\beta_k+\gamma_k$$
which expanded yields an identity of quadratic polynomials in $x$, and by identifying the coefficients we obtain the following equations for $\alpha_k,\beta_k,\gamma_k$:
\be\label{e1}
k+1+D_1=(2k+1+D_1)\, (2k+2+D_1)\, \alpha_k
\ee
\be\label{e2}
D_2=\left[ (4k+2)D_2+D_1D_2+D_2D_1\right]\alpha_k\, +\, (2k+D_1)\, \beta_k
\ee
\be\label{e3}
k-1=\left[ D_2^2-(2k+2+D_1)\right]\,\alpha_k\, +\, D_2\beta_k\, +\, \gamma_k
\ee
The system (\ref{e1})-(\ref{e2}) has a unique solution $\alpha_k,\beta_k,\gamma_k$ due to the nonresonance condition (\ref{non_res}). \qed

\subsubsection{The dominant coefficient}\label{domcoeff}
\begin{Proposition}\label{dc}
The coefficient of $x^k$ in $\bfP_k(x)$ is
\be\label{domcf}
(D_1+k+1)\, (D_1+k+2)\, \ldots\, (D_1+2k)
\ee
which is invertible under the assumptions of \S\ref{notat}.
\end{Proposition}

{\em{Proof.}} Retaining only the dominant coefficients in the representation (\ref{prodPk}), (\ref{defAj}) we have:
$$\bfP_k(x)\, =\, \left[ x(2+D_1)+x^2\partial_x\right]\,\ldots\, \left[ x(2k-2+D_1)+x^2\partial_x\right]\, \left[ x(2k+D_1)\right]\, +\, O(x^{k-1})$$
which yields
$$\bfP_k(x)=(D_1+k+1)\, (D_1+k+2)\, \ldots\, (D_1+2k)\, x^k\, +\, O(x^{k-1})$$
giving (\ref{domcf}). 

Under the assumptions and notations of (i), (ii) of \S\ref{notat}, we can assume $M_1$ diagonal, and the operator $D_1$ has the eigenvector/values:
$$D_1\mathbf{b}_{\bfm ,j}\, =\, (\mathbf{m}\cdot {\bflam}-\lambda_j )\, \mathbf{b}_{\bfm ,j}$$
Then $D_1+j$ is invertible for all $j\in\ZZ_+$ due to the assumption (ii), and therefore (\ref{domcf}) is invertible. \qed

\subsubsection{Completeness}\label{Completeness} The set $\bfP_k(x)$, $k=0,1,2,\ldots$ is complete in $\mathcal{P}_n[x]$ in the following sense:

\begin{Proposition}\label{compl}

For any $\bff=\bff(x,\bfw)\in\mathcal{P}_n[x]$ polynomial, homogeneous degree $n$ in $\bfw$ and degree $k$ in  $x$ there exist $\bfq_0,\ldots\bfq_k\in\mathcal{P}_n$ so that
\be\label{lin_comb}
\bff(x,\bfw)\, =\, \sum_{j=0}^k\, \bfP_j(x)\bfq_j(\bfw)
\ee
and the representation (\ref{lin_comb}) is unique.
\end{Proposition}

{\em{Proof.}} The decomposition (\ref{lin_comb}) follows easily by induction on $k$, relying on the fact that the dominant coefficients (\ref{domcf}) are invertible. \qed

\subsubsection{Integral inter-relation}\label{int:intrel}

Consider the following variation of the Fuchsian equation (\ref{Fuchs_eq}):
\be\label{tildeM}
\tilde{M}=M-\frac{2x}{(n-1)Q}\,I
\ee
and let $\tilde{Y}$ be a corresponding fundamental matrix: $\tilde{Y}'=\tilde{M} \tilde{Y}$. Let $\tilde{P}_{k}$ denote the polynomials  associated by (\ref{defPk}).

We have the following analogue of an integral relation for Jacobi polynomials (see \cite{Abramowitz} \S 22.13.1)
\be\label{inteq1}
{\bfP}_k(x)\bfq(\bfw)=Q(x)^{-1}Y(x)\, \int_{-1}^x\, Y(t)^{-1}\tilde{\bfP}_{k+1}(t)\bfq\left(Y(t)Y(x)^{-1}\bfw\right)\, dt
\ee
or, in differential form,
\be\label{inteq2}
\partial_x{\bfP}_k\bfq+{\rm{d}}_w({\bfP}_k\bfq) \, M\bfw-M\, \bfP_k\bfq=Q^{-1}\tilde{\bfP}_{k+1}\bfq
\ee
or, in operator notation,
\be
\left( xD_1+D_2+Q\partial_x\right)\bfP_k=\tilde{\bfP}_{k+1}
\ee

This can be easily seen due to the fact that $\tilde{\bfP}_k$  also satisfy
$$\tilde{\bfP}_{k+1}(x)\bfq(\bfw)=Q\, Y(x)\, \frac{d^{k+1}}{dt^{k+1}}\Big| _{t=x}\, \left[\, Q(t)^k\, Y^{-1}(t)\, \bfq\left(Y(t)Y^{-1}(x)\bfw\right)\,\right]$$
(which follows by a short calculation using the fact that $\bf{q}(\bfw)$ is a homogeneous polynomial in 
$\bfw$, degree $n$.)

\subsubsection{Orthogonality}\label{ortho}
The following quasi-orthogonality relations hold:
\be\label{orr}
\int_{-1}^1\, {\mathbf{P}}_j(x)\, W(x)\, \bfP_k(x)\, dx\, =\, 0\ \ {\mbox{for\ }} j<k
\ee
and
\be\label{orl}
\int_{-1}^1\, {W(x)\, \mathbf{P}}_j(x)\, \bfP_k(x)\, dx\, =\, 0\ \ {\mbox{for\ }} j>k
\ee
provided that the integrals exist.

{\em{Proof.}} For any index $k$ we have
$$[W(x)\, \bfP_k(x)]\bfq(\bfw)\, =\, Y(x)^{-1}\bfP_k(x)\bfq\left(Y(x)\bfw\right)=\frac{d^k}{dx^k} \,\left[\, Q(x)^k\, Y^{-1}(x)\, \bfq\left(Y(x)\bfw\right)\,\right]$$
and (\ref{orr}),(\ref{orl}) follow using integration by parts (which holds for matrix multiplication). \qed

\section{Acknowledgements}

The author is grateful for the warm and illuminating discussions with  Elena Berdysheva, Herb Clemens, Yuji Kodama, Doron Lubinsky, Irina Nenciu and Paul Nevai, and for the helpful and illuminating correspondence with Mourad Ismail.


\begin{thebibliography}{99}

\bibitem{Simon} B. Simon, {\em{The Analytic Theory of Matrix Orthogonal Polynomials}}, Surveys in Approximation Theory, Vol. 4, 2008, pp.1-85

\bibitem{Grunbaum} F. A. Gr\"unbaum, {\em{Matrix valued Jacobi polynomials}}, Bulletin des Sciences Mathematiques, Volume 127, Number 3, May 2003 , pp. 207-214(8)

\bibitem{RDC_Euler} R. D. Costin,   {\em{Integrability Properties of Nonlinearly Perturbed Euler Equations}}, Nonlinearity Vol. 10 No 4 pp. 905-924 (1997) 

\bibitem{Norm_Form} R. D. Costin,   {\em{Nonlinear perturbations of Fuchsian systems: correction and linearization, normal forms}}, submitted

\bibitem{RDC_MDK}  R. D. Costin, M. D. Kruskal, {\em{Nonintegrability criteria for a class of differential equations with two regular singular points}}, Nonlinearity, 16 (2003), no. 4, p.1295--1317

\bibitem{RDC_AnLin}  R. D. Costin,  {\em{Analytic linearization of nonlinear perturbations of Fuchsian systems}}, submitted


\bibitem{Abramowitz} M. Abramowitz, I. A. Stegun, {\em{Handbook of Mathematical Functions}},  National Bureau of Standards, Applied Mathematics Series - 55, June 1964










\end{thebibliography}
\end{document}